\documentclass[aop,MSNbibl,seceqn,nameyear,dvips]{arximspdf}
\usepackage{mathbh}

\doi{10.1214/11-AOP693}
\volume{40}
\issue{6}
\pubyear{2012}
\firstpage{2690}
\lastpage{2704}

\makeatletter

\newcommand{\xrightarrow}{\mathop{\hbox to 1cm{\rightarrowfill}}\limits_}

\renewcommand{\vec}[1]{\mathbf{#1}}
\newcommand{\inte}{\operatorname{Int}}

\newcommand{\schw}{\stackrel{\mathcal{D}}{\longrightarrow}}

\newtheorem{theorem}{Theorem}[section]
\newtheorem{lem}[theorem]{Lemma}
\newproclaim{rem}[theorem]{Remark}

\newcommand{\er}{\mathbb{R}}
\newcommand{\en}{\mathbb{N}}

\makeatother

\begin{document}
\begin{frontmatter}

\title{Distributions on unbounded moment spaces and random moment
sequences\thanksref{T1}}
\runtitle{Distributions on unbounded moment spaces}

\thankstext{T1}{Supported in part by the Sonderforschungsbereich
Tr/12 (Project C2, Fluctuations
and universality of invariant random matrix ensembles).}

\begin{aug}
\author[A]{\fnms{Holger} \snm{Dette}\corref{}\ead[label=e1]{holger.dette@rub.de}}
\and
\author[B]{\fnms{Jan} \snm{Nagel}\ead[label=e2]{jan.nagel@ma.tum.de}}
\runauthor{H. Dette and J. Nagel}
\affiliation{Ruhr-Universit\"at Bochum and Technische Universit\"at M\"unchen}
\address[A]{Fakult\"at f\"ur Mathematik\\
Ruhr-Universit\"at Bochum\\
44780 Bochum\\
Germany\\
\printead{e1}} 
\address[B]{Zentrum Mathematik\\
Technische Universit\"at M\"unchen\\
85747 Garching\\
Germany\\
\printead{e2}}
\end{aug}

\received{\smonth{7} \syear{2010}}
\revised{\smonth{5} \syear{2011}}

%
\begin{abstract}
In this paper we define distributions on moment spaces corresponding to
measures on the real line with an unbounded support. We identify these
distributions as limiting distributions of random moment vectors
defined on compact moment spaces and as distributions corresponding to
random spectral measures associated with the Jacobi, Laguerre and
Hermite ensemble from random matrix theory. For random vectors on the
unbounded moment spaces we prove a central limit theorem where the
centering vectors correspond to the moments of the Marchenko--Pastur
distribution and Wigner's semi-circle law.
\end{abstract}

%

\end{frontmatter}

\section{Introduction}
For a set $T \subset\mathbb{R}$, let $ \mathcal{P} (T) $ denote the
set of all probability measures on the Borel field of
$T$ with existing moments.
For a measure $\mu\in\mathcal{P} (T)$, we denote by
\[
m_k(\mu) = \int_T x^k \mu(dx);\qquad k=0,1,2,\ldots,
\]
the $k$th moment and define
%
%
\begin{equation}
\label{mominfty}
\mathcal{M} (T) = \{ \vec{m}(\mu) = (m_1(\mu),m_2(\mu),\ldots)^T
|\mu\in\mathcal{P} (T) \} \subset\mathbb{R}^{\mathbb{N}}
\end{equation}
as the set of all moment sequences. We denote by $\Pi_n$ ($n \in
\mathbb{N}$) the
canonical projection onto the first $n$ coordinates and call
$\mathcal{M}_n (T) = \Pi_n ( \mathcal{M} (T) ) \subset
\mathbb{R}^n$
the $n$th moment space. These moments
have found considerable interest in the literature; see \citet
{karstu1966}. Most authors concentrate on the ``classical'' moment space
corresponding to measures on the interval $[0,1]$; see \citet
{karsha1953}, \citet{krenud1977}, among others. \citet{chakemstu1993} equipped
it with a uniform distribution in
order to understand more fully the shape and the structure. In
particular, these authors proved asymptotic normality of an appropriately
standardized version of a projection $\Pi_k (\vec{m}_n)$ of a
uniformly distributed vector $\vec{m}_n$ on $\mathcal{M}_n([0,1])$.
\citet{gamloz2004} considered large deviation principles for
random moment sequences on the space $\mathcal{M}_n([0,1])$, while
\citet
{lozada2005}
investigated similar problems for moment spaces corresponding to more
general functions defined on a bounded set. More recently, \citet{gamrou2010}
discussed random spectral measures related to moment spaces of
measures on the interval $[0,1]$ and
moment spaces related to measures defined on the unit circle.
Distributions of random moments induced by measures on $l_p$ balls were
investigated by \citet{bartheetal2010}.

The present paper is devoted to the problem of defining probability
distributions on unbounded moment spaces. We will investigate these
distributions from several perspectives. We introduce a class of
general distributions on the moment space corresponding to measures
defined on a compact
interval. By a limiting argument we will derive canonical distributions
on the moment spaces corresponding to measures on the unbounded intervals
$T=[0,\infty)$ and $\mathbb{R}$, respectively. We also show that these
distributions appear naturally in the study of random spectral measures
of the classical Jacobi, Laguerre and Gaussian ensemble.
Finally we consider random moment sequences distributed according to
the new probability
distributions on the unbounded moment spaces. In particular, we prove
weak convergence of a centered random moment vector, where the centering
vector corresponds to the moments of the Marchenko--Pastur law (in the
case of the moment space $\mathcal{M}([0, \infty))$) and to
the semi-circle law [for the moment space $\mathcal{M} (\mathbb{R})$].
These measures play a very important role in free probability and
quantum probability; see the books of \citet{hiaipetz2000} and \citet
{horaobata2007}.

\section{Distributions on unbounded moment spaces}\label{sec2}


In the present section we will
define a class of more general distributions on the $n$th moment space
$\mathcal{M}_n([a,b])$ as considered by \citet{chakemstu1993}. The
motivation for considering this class is twofold. On one hand, we want to
introduce distributions on the moment space $ \mathcal{M}_{k}([a,b])$,
which are different from the uniform distribution. On the other hand,
we want to identify distributions on unbounded moment spaces as limits
of distributions on $ \mathcal{M}_{k}([a,b])$, when $b-a \to\infty$.

Let $\vec{m}_{k-1}=(m_1,\ldots,m_{k-1})^T \in\mathcal{M}_{k-1}([a,b])$
be a
given vector of moments of a probability measure on the interval
$[a,b]$, then these
first $k-1$ moments impose bounds on the $k$th moment $m_{k}$ such that
the moment vector $\vec{m}_{k}=(m_1,\ldots,m_{k-1},m_{k})^T$ is an
element of the $k$th moment space $ \mathcal{M}_{k}([a,b])$. More
precisely, define for $\vec{m}_{k-1} \in\mathcal{M}_{k-1}([a,b]) $
\begin{eqnarray*}
m_{k}^- &=& \min\biggl\{ m_{k}(\mu) \Big| \mu\in\mathcal{P}([a,b])
\mbox{ with } \int_a^b t^i \,d\mu(t) = m_i \mbox{ for } i = 1,\ldots,k-1
\biggr\}, \\
m_{k}^+ &=& \max\biggl\{ m_{k}(\mu) \Big| \mu\in\mathcal
{P}([a,b]) \mbox{ with } \int_a^b t^i \,d\mu(t) = m_i \mbox{ for } i
= 1,\ldots,k-1 \biggr\},
\end{eqnarray*}
then it follows that $ \vec{m}_k
= (m_1,\ldots,m_{k})^T \in\inte\mathcal{M}_{k}([a,b]) $ if and only if
$ m_{k}^- <m_k < m_{k}^+$, where $\inte C$ denotes the interior of a
set $C \subset\mathbb{R}^k$. Consequently, we define for a point
$\vec{m}_k \in\inte\mathcal{M}_{k}([a,b])$ the canonical moment of
order $l=1,\ldots,k$ as
%
%
\begin{equation}\label{canmom}
p_l =p_l(\vec{m}_k) = \frac
{m_{l}-m_{l}^-}{m_l^+-m_l^-};\qquad l=1,\ldots,k.
\end{equation}
Note that for $\vec{m}_k\in\inte\mathcal{M}_k ([a,b])$ we have
$p_l\in(0,1); l=1,\ldots,k$; and that $p_k$ describes
the relative position
of the moment $m_k$ in the set of all possible $k$th moments with
fixed moments $m_1,\ldots,m_{k-1}$. For this reason, the canonical
moments do not depend on the interval $[a,b]$, that is they are
invariant under linear transformations of the measure; see \citet
{dettstud1997}. Moreover, definition (\ref{canmom})
defines a one-to-one mapping,
%
%
\begin{equation}\label{phimap}
\varphi_n \dvtx\vec{m}_n \mapsto\vec{p}_n =(p_1,\ldots,p_n)^T
\end{equation}
from the interior of the moment space $\mathcal{M}_n([a,b]) $ onto the
open cube $(0,1)^n$.
It can be shown that for a point
$(m_1,\ldots,m_{2n-1}) \in\inte\mathcal{M}_{2n-1} ([a,b])$ the
canonical moments appear in the three-term recurrence relation
%
%
\begin{equation}
\label{polrek}\qquad
x P_k(x) = P_{k+1}(x) + b_{k+1} P_k(x) + a_k
P_{k-1}(x),\qquad
k=1,\ldots,n-1,
\end{equation}
[$P_0(x) = 1,P_1(x) = x-b_1$] of the monic
orthogonal polynomials 
associated with the vector $(m_1,\ldots,m_{2n-1})$; see \citet
{chihara1978}. These polynomials are orthogonal
with respect to every measure with first moments $m_1,\ldots,m_{2n-1}$,
and the recursion coefficients in (\ref{polrek}) are given by
%
%
\begin{eqnarray}
\label{chain1}\qquad
b_{k+1} &=& a + (b-a) \bigl((1-p_{2k-1}) p_{2k} + (1-p_{2k})
p_{2k+1}\bigr);\nonumber\\[-8pt]\\[-8pt]
&&\eqntext{k=0,\ldots,n-1,} \\
\label{chain2}
a_k &=& (b-a)^2 (1-p_{2k-2}) p_{2k-1} (1-p_{2k-1})
p_{2k};\qquad k=1,\ldots,n-1,
\end{eqnarray}
where we put $p_{-1}=p_0=0$
(note that $a_k >0; k=1,\ldots,n)$. In the case $T=[0,\infty)$
the upper bound $m_k^+$ does not exist, but we can still define for a
point $\vec{m}_{k-1} \in\inte\mathcal{M}_{k-1} ([0,\infty))$ the
lower bound
\begin{eqnarray*}\hspace*{3pt}
m_{k}^- = \min\biggl\{ m_{k}(\mu) \Big| \mu\in\mathcal
{P}([0,\infty
)) \mbox{ with } \int_0^{\infty} t^i \,d\mu(t) = m_i \mbox{ for } i =
1,\ldots
,k-1 \biggr\},
\end{eqnarray*}
where $\vec{m}_k=(m_1,\ldots,m_{k})^T
\in
\inte\mathcal{M}_{k}([0,\infty)) $ if and only if $ m_{k} > m_{k}^-$. In
this case, the analogs of the canonical moments are defined by the
quantities
%
%
\begin{equation}\label{defz}
z_l = z_l(\vec{m}_n) =
\frac{m_{l}-m_{l}^-}{m_{l-1}-m_{l-1}^-},\qquad l=1,\ldots,k
\end{equation}
(with
$m_0^-=0$) and related to the coefficients in the three-term recurrence relation
(\ref{polrek}) for the monic orthogonal polynomials by
%
%
\begin{equation}\label{zerlegungan}
a_k = z_{2k-1} z_{2k},\qquad
b_k = z_{2k-2} + z_{2k-1} .
\end{equation}
Note that (\ref{defz}) defines a one-to-one mapping
%
%
\begin{equation} \label{defpsi}
\psi_n \dvtx\vec{m}_n \mapsto\vec{z}_n = (z_1,\ldots,z_n)^T
\end{equation}
from the interior of the moment space $\mathcal{M}_n([0,\infty)) $ onto
$(\er^+)^n$. Finally, in the case $T= \mathbb{R}$, neither $m_{k}^- $
nor $m_{k}^+$ can be defined. Nevertheless, there exists also a
one-to-one mapping
%
%
\begin{equation}\label{defxi}
\xi_n \dvtx\vec{m}_{2n-1} \mapsto
(b_1,a_1,\ldots,a_{n-1},b_n)^T
\end{equation}
from the interior of the $(2n-1)$th moment space $\inte\mathcal
{M}_{2n-1}(\er)$ onto the space $( \er\times\er^+ )^{n-1}
\times\er$ of coefficients in the three-term recurrence relation~(\ref{polrek}), which
can be considered as the analog of (\ref{defpsi}) and is defined~by
%
%
\begin{eqnarray}
\label{wall1}
\int_\mathbb{R} x^k P_k(x) \,d\mu(x) &=& a_1\cdots a_k,\qquad k=1,\ldots
,n-1 , \\
\label{wall2}
\int_\mathbb{R} x^{k+1} P_k(x) \,d\mu(x) &=& a_1\cdots a_k (b_1 + \cdots
+b_{k+1}),\qquad k=0,\ldots,n-1;\hspace*{-32pt}
\end{eqnarray}
see, for example,
\citet{wall1948}. 

In the following sections we will use the canonical moments and
corresponding quantities on the interval $[0,\infty)$ and the
real line for the definition of distributions on the corresponding
moment spaces.
{The basic idea is to define a general class of distributions on the
moment space $\mathcal{M}_n([a,b])$
and to consider the limit as $b-a\to\infty$. To be precise,
let for $k \geq1$
$ \tilde{f}_k\dvtx(0,1)
\longrightarrow\er
$ be a nonnegative integrable function with $\int_0^1\tilde
{f}_k(x)\,dx >0$,}
then a distribution on the interior of the moment space $ \mathcal
{M}_{n}([a,b])$ is given by
%
%
\begin{equation}\label{defallg}
f_n(\vec{m}_n) = c_n
\prod_{k=1}^n \tilde{f}_k(p_k(\vec{m}_n)) \mathbh{1}_{ \{
m_k^-<m_k<m_k^+ \} },
\end{equation}
where $p_k(\vec{m}_n)$ is the $k$th canonical moment defined in (\ref
{canmom}) and $c_{n}$
a normalization constant. 
Our first theorem gives the distribution of the canonical moments
corresponding to the random vector $\vec{m}_n$ with density $ f_n$
defined in (\ref{defallg}).
\begin{theorem} \label{vertkanmom}
Suppose that $\vec{m}_n $ is a random vector on the moment space
$\mathcal{M}_n([a,b])$ with density $ f_n$ defined in
(\ref{defallg}). Then\vadjust{\goodbreak} the canonical moments $p_1(\vec{m}_n),\ldots
,p_n(\vec{m}_n)$ are independent and $p_k(\vec{m}_n)$ has a density
proportional to
$ \tilde{f}_k(x) (x-x^2)^{n-k} \mathbh{1}_{ \{ 0<x<1 \} }$
$(1 \leq k \leq n)$.
\end{theorem}
\begin{pf}
It follows from Theorem 1.4.9 and equation (1.3.6) in
\citet{dettstud1997} that the Jacobian
determinant of the mapping $\varphi_n$ defined in (\ref{phimap}) is
given by
%
%
\begin{eqnarray} \label{jacobiab}
\biggl| \frac{\partial\varphi_n}{\partial\vec{m}_n} \biggr|
&=& \prod_{k=1}^n \frac{\partial p_k(\vec{m}_n)}{\partial m_k}
= \prod_{k=1}^n (m_k^+-m_k^-)^{-1}
\nonumber\\[-8pt]\\[-8pt]
&=& (b-a)^{-n(n+1)/2} \prod_{k=1}^n \bigl( p_k(1-p_k)\bigr)^{-(n-k)};\nonumber
\end{eqnarray}
considering the product structure of $f_n$, this gives the asserted
distribution.
\end{pf}

For the construction of distributions on the unbounded moment
space $\mathcal{M}_n([0$, $\infty))$, a special case will be of particular
interest; that is, $\tilde{f}_k(x) = x^{\gamma_k} (1-x)^{\delta_k}$,
where $\gamma=
(\gamma_k)_{k \in\mathbb{N}}$, $\delta= (\delta_k)_{k \in\mathbb
{N}}$ are sequences of real parameters, such that $\gamma_k,\delta_k
>0$ for all $k \geq1$.
In this case the density on
the moment space $ \mathcal{M}_{n}([a,b])$ is given by
%
%
\begin{equation}
\label{momdichteab}
f_n^{(\gamma, \delta)}(\vec{m}_n) = c_n^{[a,b]}
\prod_{k=1}^n \biggl( \frac{m_k - m_k^-}{m_k^+ - m_k^-}
\biggr)^{\gamma
_k} \biggl( \frac{m_k^+ - m_k}{m_k^+ - m_k^-} \biggr)^{\delta_k}
\mathbh{1}_{\{ m_k^-<m_k<m_k^+ \} } ,\hspace*{-35pt}
\end{equation}
where
%
%
\begin{equation} \label{normconst}
c_n^{[a,b]} = \Biggl[ (b-a)^{n(n+1)/2} \prod_{k=1}^n \int^1_0
x^{n-k+\gamma_k} (1-x)^{n-k+\delta_k} \,dx \Biggr]^{-1}
\end{equation}
is the normalizing constant. The choice of density (\ref{momdichteab})
is motivated by results of \citet{dettstud1995} who showed that the
empirical distribution of the (appropriately normalized) roots of the
Jacobi polynomials $P_k^{(\gamma_k, \delta_k)} (x)$ converges weakly to
a distribution with
unbounded support if $\gamma_k \to\infty$ or $\delta_k \to\infty$.
It follows from Theorem~\ref{vertkanmom} that for the density
$f_n^{(\gamma,\delta)}$, the canonical moment $p_k$ has a Beta
distribution $\operatorname{Beta}(\gamma_k +n-k+1,\delta_k+n-k+1)$.
In the following we use
densities of form (\ref{momdichteab}) to construct a distribution on
the unbounded moment space
$ \mathcal{M}_n ([0,\infty))$.
\begin{theorem} \label{momkonv1}
Let $f_n^{(\gamma^{(d)} ,\delta^{(d)})}$ denote the density defined in
(\ref{momdichteab}) on the moment space $\mathcal{M}_n([0,d])$
corresponding to the probability measures on the interval $[0,d]$, where
the parameter sequences $\gamma^{(d)} =(\gamma^{(d)}_k)_{k \in
\mathbb
{N}}, \delta^{(d)} = (\delta_k^{(d)})_{k \in\mathbb{N}}$ depend on
length $d$ and satisfy for all $k \geq1$
$ \gamma^{(d)}_k \rightarrow\gamma_k > -1$ and ${\delta^{(d)}_k}/{d}
\rightarrow\delta_k \in\er^+
$
as $d \rightarrow\infty$. Then for $d \to\infty$ the density
$f_n^{(\gamma^{(d)} ,\delta^{(d)})}$ converges point-wise to the
function
%
%
\begin{equation}\label{momdichte0infty}\qquad
g_n^{(\gamma,\delta)}(\vec{m}_n)
= c_n^{[0,\infty)} \prod_{k=1}^n z_k(\vec{m}_n)^{\gamma_k} \exp
(-\delta_k z_k(\vec{m}_n)) \mathbh{1}_{ \{ z_k(\vec{m}_n) >0 \} },
\end{equation}
where $ z_k(\vec{m}_n)$ is given in (\ref{defz}) and
$c_n^{[0,\infty)} = \prod_{k=1}^n ( {\delta_k^{\gamma_k +n-k+1}})
/\Gamma(\gamma_k +n- k+1) $.
Moreover, $g_n^{(\gamma,\delta)}$ defines a density on the unbounded
moment space $\mathcal{M}_n ([0,\infty))$.
\end{theorem}
\begin{pf}
The fact that $g_n^{(\gamma,\delta)}$ is a density is
obvious from the transformation in the proof of Theorem~\ref{vertzeta}
below, we prove here only the convergence. For a fixed point $\vec{m}_n
\in\mathcal{M}_n ([0,\infty))$,
there exists a $d_0 \in\mathbb{N}$ with $\vec{m}_n \in\mathcal{M}_n
([0,d])$ for all $d \geq d_0$. Let
$\vec{p}_n(\vec{m}_n)$ denote the vector of canonical moments
corresponding to the vector $\vec{m}_n$ in the moment space $\mathcal
{M}_n([0,d])$.
We will show at the end of this proof that
%
%
\begin{equation} \label{expan1}
p_k(\vec{m}_n) = \frac{z_k(\vec{m}_n)}{d} \bigl(1+ o(1)\bigr),
\qquad k=1,\ldots,n,
\end{equation}
where the quantities $z_k(\vec{m}_n)$ are defined in (\ref{defz}).
Observing this representation and definition (\ref{momdichteab}), it
follows for $d\to\infty$
\begin{eqnarray*}f_n^{(\gamma^{(d)}
,\delta^{(d)})}(\vec{m}_n) &=& c_n^{[0,d]} \prod_{k=1}^n \biggl(
\frac
{z_k(\vec{m}_n)}{d} \biggr)^{\gamma^{(d)}_k}
\biggl( 1- \frac{z_k(\vec{m}_n)}{d} \biggr)^{\delta_k^{(d)}}
\bigl(1+o(1)\bigr) \\
&=& d^{-(\gamma^{(d)}_1+\cdots+\gamma^{(d)}_n)} c_n^{[0,d]} \prod
_{k=1}^n z_k(\vec{m}_n)^{\gamma_k} \exp(-\delta_k z_k(\vec{m}_n))
\bigl(1+o(1)\bigr) .
\end{eqnarray*}
Finally, we obtain from (\ref{normconst}) for the normalizing constant
by a straightforward calculation
$
d^{-(\gamma^{(d)}_1+\cdots+\gamma^{(d)}_n)} c_n^{[0,d]}
= c_n^{[0,\infty)} (1+o(1)) ,
$
which proves the assertion of the theorem. For the remaining proof of
the representation (\ref{expan1}), let $\mu$ be a measure on the
interval $[0,d]$
with first moments given by $\vec{m}_n$ and let $\nu$ denote the
measure on the interval $[0,1]$ obtained from $\mu$ by the linear
transformation $x \mapsto x/d$.
We write $p_k (\mu)$ for $p_k(\vec{m}_n)$ and $z_k(\mu)$ for
$z_k(\vec
{m}_n)$. Invariance of the canonical moments under linear
transformations yields $p_k(\mu) = p_k(\nu)$. The recursion variables
of the measure $\nu$ can be decomposed as
%
%
\begin{equation} \label{zetakanmom}
z_k(\nu) = \bigl(1-p_{k-1}(\nu)\bigr) p_k(\nu) .
\end{equation}
A comparison of the continued fraction expansion of the Stieltjes
transform of $\mu$ and of $\nu$ [see Theorem 3.3.3 in \citet
{dettstud1997}] yields
$dz_k(\nu) = z_k(\mu)$. With (\ref{zetakanmom}) we obtain
\[
p_k(\mu) = p_k(\nu) = \frac{z_k(\nu)}{ 1- p_{k-1}(\nu)} = \frac{1}{d}
\frac{z_k(\mu)}{ 1- p_{k-1}(\mu)}
\]
for $k >1$. The first canonical moment is given by
$
p_1(\vec{m}_n) = {(m_1-m_1^-)}/(m_1^+ - m_1^-) = {m_1}/{d} =
{z_1(\vec
{m}_n)}/{d},
$
and (\ref{expan1}) follows by an induction argument.
\end{pf}

The following theorem is essential for the asymptotic
investigations in Section~\ref{sec3} and gives the distribution of the the
vector $\vec{z}_n=(z_1,\ldots,z_n)^T$
corresponding to a random point in $\mathcal{M}_n ([0,\infty))$.
%
\begin{theorem} \label{vertzeta}
Let $\vec{m}_n \in\mathcal{M}_n([0,\infty))$ be governed by a law with
density $g_n^{(\gamma,\delta)}$. Then the recursion variables
$\vec{z}_n = \psi_n (\vec{m}_n)$ defined by (\ref{defz}) are
independent and gamma distributed, that is,
\[
z_k \sim \operatorname{Gamma}(\gamma_k+n-k+1,\delta_k),\qquad k=1,\ldots,n.
\]
\end{theorem}
\begin{pf}
By its definition in (\ref{defz}), the random variable
$z_k$ depends only on the moment $m_1,\ldots,m_{k}$;
therefore the Jacobi matrix of the mapping $\psi_n$ is a lower triangular
matrix. We obtain for the Jacobian determinant
\[
\biggl| \frac{\partial\vec{m}_n}{\partial\vec{z}_n} \biggr| =
\prod
_{k=1}^n \biggl| \frac{\partial m_k}{\partial z_k} \biggr|
= \prod_{k=1}^n (m_{k-1} - m_{k-1}^-) = \prod_{k=2}^n z_1 \cdots
z_{k-1}= \prod_{k=1}^n
z_k^{n-k},
\]
where the third identity follows from the definition of the $z_i$ in
(\ref{defz}).
Considering the second representation of the density in (\ref
{momdichte0infty}), this yields the claimed distribution.
\end{pf}

We conclude this section with a discussion of distributions on the
moments space corresponding to measures on $\er$. For the sake of
brevity we
restrict ourselves to moment spaces of odd dimension, that is,
$\mathcal{M}_{2n-1} (\er)$.
To derive a class of distributions on $\mathcal{M}_{2n-1} (\er)$
we consider the moment space $\mathcal{M}_{2n-1}([-s,s])$ with $s \to
\infty$ and a density of the form (\ref{momdichteab}) with parameters
varying with $s$.
The proof of the following result is similar to the proof of Theorem
\ref{momkonv1} and therefore omitted.
\begin{theorem} \label{konvdichte3}
Denote by ${f}_{2n-1}^{(\gamma^{(s)}, \delta^{(s)} )}$ the density
defined in (\ref{momdichteab}) on the moment space $\mathcal
{M}_{2n-1}([-s,s])$, where the
parameters satisfy
\begin{eqnarray*}
\gamma_{2k-1}^{(s)} &=& \delta_{2k-1} s^2 + o(1),\qquad \delta
_{2k-1}^{s} = \delta_{2k-1} s^2 + o(1), \\
\gamma_{2k}^{(s)} &=& \gamma_k + o(1),\qquad \delta_{2k}^{(s)} =
\delta
_{2k} s^2 + o(s^2)
\end{eqnarray*}
with $\gamma_k > -1, \delta_k> 0$. Then
$f_{2n-1}^{(\gamma^{(s)}, \delta^{(s)} )}$ converges point-wise to
the function
%
%
\begin{eqnarray} \label{momdichteer}
&&h_{2n-1}^{(\gamma,\delta)}(\vec{m}_{2n-1})\nonumber\\
&&\qquad= \prod_{k=1}^n \sqrt
{\frac{\delta_{2k-1}}{\pi}} \exp( -\delta_{2k-1} b_k^2 (\vec
{m}_{2n-1}) ) \nonumber\\[-8pt]\\[-8pt]
&&\qquad\quad{} \times\prod_{k=1}^{n-1} \frac{{\delta_{2k}}^{\gamma_k
+2n-2k}}{\Gamma(\gamma_k +2n-2k)} a_k^{\gamma_k}
(\vec{m}_{2n-1}) \nonumber\\
&&\qquad\quad\hspace*{28.3pt}{}\times\exp( -\delta_{2k} a_k (\vec{m}_{2n-1})
)
\mathbh{1}_{ \{ a_k (\vec{m}_{2n-1}) >0 \} } . \nonumber
\end{eqnarray}
Moreover,
the function $h_{2n-1}^{(\gamma,\delta)}$ defines a density on the
moment space $\mathcal{M}_{2n-1}(\er)$.
\end{theorem}

The following result is the
analog of Theorem~\ref{vertzeta}. The proof follows by similar
arguments, where the Jacobian can be obtained from equations (\ref
{wall1}) and (\ref{wall2}).
\begin{theorem} \label{vertanbn}
Let $\vec{m}_{2n-1} \in\mathcal{M}_{2n-1}(\er)$ be a random vector
with density $h_{2n-1}^{(\gamma,\delta)}$ defined in (\ref{momdichteer}).
Then the random recursion coefficients $(b_1,a_1,\ldots,\break a_{n-1}$,
$b_n)^T =
\xi_{2n-1} ( \vec{m}_{2n-1})$ in the recurrence relation (\ref{polrek})
for the orthogonal polynomials associated with $\vec{m}_{2n-1}$ are
independent and
\[
b_k \sim\mathcal{N} \biggl(0,\frac{1}{2\delta_{2k-1}}\biggr),
\qquad a_k \sim \operatorname{Gamma}
(\gamma_k +2n-2k,\delta_{2k}).
\]
\end{theorem}

\begin{rem} \label{ranspect}
It is notable that the introduced distributions on the moment space
appear naturally as
distributions of moment vectors corresponding to random spectral
measures which were recently
discussed by Gamboa and Rouault
(\citeyear{gamrou2009,gamrou2010}). To be precise, let
$\vec{w}=(w_1,\ldots,w_n)$ and $\bolds{\lambda} =(\lambda_1,\ldots,
\lambda_n)$ denote independent random variables.
Assume that $\vec{w}$ follows a Dirichlet distribution with density
%
%
\begin{equation}
\label{dirdensity}
\frac{\Gamma({n\beta}/{2} )}{\Gamma({\beta}/{2})^n}
w_1^{\beta/2-1} \cdots w_n^{\beta/2-1} \mathbh{1}_{\{ w_1,\ldots,w_n
>0, \sum_{i=1}^n w_i =1 \} }
\end{equation}
and that the density of $\bolds{\lambda}$ is given by the Jacobi-ensemble
%
%
\begin{equation} \label{jacew}
f_J(\lambda) = c_J^{\gamma_0, \delta_0} | \Delta(\lambda)|^{\beta}
\prod_{i=1}^n \lambda_i^{\gamma_0} (1-\lambda_i)^{\delta_0}
\mathbh
{1}_{\{ 0<\lambda_i<1 \} } ,
\end{equation}
where $\gamma_0, \delta_0 > - \frac{1}{2}$; see
\citet{kilnen2004}. We define a random measure by $ \mu=
\sum_{i=1}^n w_i \delta_{\lambda_i}$, and it follows from Theorem 2.2
and Proposition 5.3 in \citet {kilnen2004} (applied to the
interval $[-1,1]$) that $\mu$ is the spectral measure of the random
tridiagonal matrix
\[
J_n =
\pmatrix{
d_1 & c_1 & & \cr
c_1 & d_2 & \ddots& \cr
& \ddots& \ddots& c_{n-1} \cr
& & c_{n-1} & d_n},
\]
where $
d_k = p_{2k-2} (1-p_{2k-3}) + p_{2k-1}(1-p_{2k-2})$ and
\[
c_k = \sqrt{
p_{2k-1}(1-p_{2k-2})p_{2k}(1-p_{2k-1})}\vadjust{\goodbreak}
\]
and $p_{-1}=p_0=0$ and $p_1,\ldots,p_{2n-1}$ are independent random
variables distributed as
\[
p_k \sim
\cases{\displaystyle \operatorname{Beta}\biggl( \frac{2n-k}{4} \beta, \frac{2n-k-2}{4}
\beta+\gamma_0+\delta_0+2 \biggr), &\quad $k$ even, \vspace*{2pt}\cr
\displaystyle \operatorname{Beta}\biggl( \frac{2n-k-1}{4} \beta+\gamma_0+1, \frac{2n-k-1}{4}
\beta+\delta_0+1 \biggr), &\quad $k$ odd.}
\]
%
This spectral measure is the unique measure with $\langle e_1,J_n^k e_1
\rangle= m_k(\mu)$ for all~$k$.
The tridiagonal matrix $J_n$ defines monic polynomials $P_1(x),\ldots
,P_n(x)$ via a recursion (\ref{polrek}) with recursion coefficients
$b_k = d_k $ ($1 \leq k\leq n$), $a_k = c_k^2$ ($1 \leq k \leq n-1$).
Indeed, the polynomial $P_k(x)$ is the characteristic polynomial of the
upper left $(k \times k)$-subblock of the matrix $J_n$ and the $k$th
orthogonal polynomial with respect do the measure $d\mu$.
Therefore, we obtain
from (\ref{chain1}) and (\ref{chain2}) 
that
$\vec{p}_{2n-1} =(p_1,\ldots,p_{2n-1})^T$ is exactly the vector of
canonical moments of the spectral measure $\mu$, and by definition,
their joint density is given by
\begin{eqnarray*}
f_{\vec{p}} (\vec{p}_{2n-1})
&=& c \prod_{k=1}^n p_{2k-1}^{(n-k)\beta/2+\gamma_0}
(1-p_{2k-1})^{(n-k)\beta/2+\delta_0}\\
&&{}\times\prod_{k=1}^{n-1} p_{2k}^{(n-k)\beta/2-1} (1-p_{2k})^{(n-k-1)\beta
/2+\gamma_0+\delta_0+1}.
\end{eqnarray*}
Since the eigenvalues of the matrix $J_n$ are contained in the interval $(0,1)$,
the moments $\vec{m}_{2n-1}(\mu) = \varphi_{2n-1}^{-1}(\vec{p}_{2n-1})$
of the spectral measure are in the interior of the
moment space $\mathcal{M}_{2n-1}([0,1])$. The Jacobian of the
transformation $\varphi_{2n-1}^{-1}$ is given by
$\prod_{k=1}^n (p_k(1-p_k))^{2n-1-k}$,
which gives for the density of the random moment vector $\vec
{m}_{2n-1}(\mu)$
\begin{eqnarray*}
&&f_{\vec{m}} (\vec{m}_{2n-1})\\
&&\qquad= c \prod_{k=1}^n p_{2k-1}(\vec{m}_{2n-1})^{(\beta/2-2)(n-k)+\gamma_0}
\bigl(1-p_{2k-1}(\vec{m}_{2n-1})\bigr)^{(\beta/2-2)(n-k)+\delta_0} \\
&&\qquad\quad{} \times\prod_{k=1}^{n-1} p_{2k}(\vec{m}_{2n-1})^{(\beta/2-2)(n-k)}
\bigl(1-p_{2k}(\vec{m}_{2n-1})\bigr)^{(\beta/2-2)(n-k-1)+\gamma_0+\delta_0}.
\end{eqnarray*}
This is a density as in (\ref{momdichteab}) with parameters
$\gamma_{2k-1} = (\frac{\beta}{2} -2)(n-k) +\gamma_0 $, $\delta_{2k-1}
= (\frac{\beta}{2} -2)(n-k) +\delta_0$,
($1\leq k \leq n$) and $\gamma_{2k} = (\frac{\beta}{2} -2)(n-k)$,
$\delta_{2k} = (\frac{\beta}{2} -2)(n-k-1) +\gamma_0+ \delta_0$
($1\leq k \leq n-1$).
We finally note that densities of the form (\ref{momdichte0infty}) and
(\ref{momdichteer}) on the moment space $\mathcal{M}_n([0,\infty)) $
and $ \mathcal{M}_{2n-1}(\er) $
are obtained
by replacing the Jacobi ensemble by the Laguerre and Hermite ensemble,
respectively [see \citet{nagdet2010}
for details].
\end{rem}

\section{Weak convergence of random moments}\label{sec3}

In this section we study the probabilistic properties of random vectors
on the moment spaces $\mathcal{M}_n([a,b])$, $\mathcal{M}_n([0,\infty))$
and $\mathcal{M}_n(\mathbb{R})$ distributed according to the measures
introduced in
Section~\ref{sec2}. We begin with random moments defined on the moment space
corresponding to probability measures on a compact interval.
\citet{chakemstu1993} and \citet{gamloz2004} investigated the uniform
distribution on $\mathcal{M}_n([a,b])$, and we first demonstrate
that weak convergence of random moment vectors can be established for a
rather broad class of distributions on $\mathcal{M}_n ([a,b])$.
An important role in the discussion of moment spaces corresponding to
probability
measures with bounded support $[a,b]$ plays the arcsine distribution
$\nu$ with density
$
d\nu(x) = {1}/{\pi\sqrt{ (x-a)(b-x)} } \mathbh{1}_{\{ a< x <b \} } \,dx.
$
The canonical moments of the arcsine distribution are given by $1/2$
[see \citet{dettstud1997}], and therefore its sequence of moments could
be considered as the
center of the moment space $\mathcal{M}([a,b])$.
The following statements establish the asymptotic properties of the
(random) canonical moments corresponding to
distributions on the moment space $\mathcal{M}_n([a,b])$ defined in
(\ref{defallg}).
Throughout\vspace*{1pt} this paper the symbol $\schw$ stands for weak convergence.
\begin{theorem} \label{allgthm1}
Suppose that the distribution of the random moment vector $\vec{m}_n
\in\mathcal{M}_n([a,b])$ is absolute continuous with density $ f_n$
defined in
(\ref{defallg}),
where the point $\frac{1}{2}$ is in the support of the measure with
density proportional to $ \tilde{f}_k$
and denote by $p_k^{(n)}$ the $k$th canonical moment of $\vec{m}_n$
$(k=1,\ldots,n)$:
\begin{longlist}[(a)]
\item[(a)] If $n \to\infty$, then almost surely 
$
p_k^{(n)} \xrightarrow{n \rightarrow\infty} \frac{1}{2}.
$
\item[(b)]
If additionally the function $\tilde{f}_k$ in the density (\ref
{defallg}) is bounded, continuous at $\frac{1}{2}$ and positive, then
the $k$th canonical moment corresponding to~$\vec{m}_n$ satisfies
$
\sqrt{8n} (p_k^{(n)} -\frac{1}{2} ) \xrightarrow
{n \rightarrow\infty}^{\mathcal{D}} \mathcal{N}(0,I).
$
\end{longlist}
\end{theorem}
%
%
\begin{pf}
For notational convenience, we consider only the case
$[a,b]=[0,1]$. The proof of the first assertion is a consequence of the
Laplace method [see, e.g., \citet{polyaszegoe1998},\vspace*{2pt} page 96]: The
canonical moment has density proportional to $\tilde{f}_k(x)
(x-x^2)^{n-k}$, which concentrates exponentially fast on\vspace*{1pt} any
neighborhood of $\frac{1}{2}$. The almost sure convergence follows
then with an application of the Borel--Cantelli lemma.

For a proof of part (b), we apply a technique similar to stable
convergence; see \citet{aldeag1978} and the references therein. From
\citet{chakemstu1993} we know that for $\tilde{f}_k \equiv1$ the
convergence holds, that is, $E[I_A(\tilde{p}_k)] \rightarrow E[I_A(p)]$
for $\tilde{p}_k=\sqrt{8n}(p_k^{(n)}-\frac{1}{2})$ the normalized
canonical moment, $p$ standard normal distributed and $A$ of the form
$(-\infty,a]$. Then for any $\tilde{f}_k$ satisfying the assumptions
of part
(b), the convergence $E[I_A(\tilde{p}_k)\tilde{f}_k(\frac
{1}{2}+\frac
{1}{\sqrt{8n}}\tilde{p}_k )] \rightarrow E[I_A(p)\tilde{f}_k(\frac
{1}{2})]$ holds. This implies the convergence of $\tilde{p}_k$ if the
density of $p_k^{(n)}$ is multiplied by $\tilde{f}_k$.
\end{pf}

\citet{chakemstu1993}
showed weak convergence of the vector $\vec{m}^{(n)}_k$ of the first
$k$ components of a uniformly distributed moment vector
$\vec{m}_n=(m_1,\ldots,m_n)$ on the moment space $\mathcal
{M}_n([0,1])$ (i.e., $f_k \equiv1$), that is,
%
%
\begin{equation} \label{chakemsturesult}
\sqrt{8n}A^{-1}\bigl(\vec{m}_k^{(n)} -\vec{m}_k(\nu)\bigr)
\xrightarrow{n\longrightarrow\infty}^{\mathcal{D}} \mathcal{N}(0,I_k),
\end{equation}
where $\vec{m}_k(\nu)$ denotes the vector of
the first $k$ moments of the arcsine distribution and $A$ is a $k
\times k$ lower triangular matrix with entries
\[
a_{i,j} = 2^{-2i+2} \pmatrix{2i\cr i-j},\qquad j \leq i .
\]
By part (b) of Theorem~\ref{allgthm1}, it is easy to see that the
weak convergence in~(\ref{chakemsturesult}) holds for the more general
densities $f_n$ on $\mathcal{M}_n([0,1])$.
By the arguments in the proof of Theorem~\ref{allgthm1} it is also
apparent that no condition on the density $f_n$ can be dropped, without
adding more specific restrictions. In particular, we need $\tilde{f}_k$
to be independent of $n$, which implies the product structure of $f_n$.
In this sense, the densities as in (\ref{defallg}) constitute the
largest class of densities on the moment space such that the
convergence (\ref{chakemsturesult}) holds.

We conclude this paper with a discussion of corresponding results for
distributions on the noncompact moment spaces
$\mathcal{M}_n([0,\infty) )$ and $\mathcal{M}_n(\er)$. In this case the
analogs of the arcsine distribution
in this context are the Marchenko--Pastur and Wigner's semicircle
distribution defined by
%
%
\begin{equation}
\label{marchpast} d\eta(x) = \frac{\sqrt{x(4-x)}}{2\pi x} \mathbh{1}_{
\{ 0 <x<4 \} } \,dx,\qquad d\rho(x) = \frac{1}{2\pi} \sqrt{4-x^2}
\mathbh{1}_{ \{ -2<x<2 \} } \,dx,\hspace*{-35pt}
\end{equation}
respectively; see \citet{hiaipetz2000} or \citet{nicspe2006}. The moments
of the Marchenko--Pastur law $\eta$ are the Catalan numbers $c_n$
defined by
$
m_n(\eta) = c_n = \frac{1}{n+1} {2n\choose n}$ ($ n \in\mathbb{N}$),
and the moments of the semicircle law $\rho$ are given by $m_{2n}(\rho
)=c_n, m_{2n-1}(\rho)=0$.
Our next result establishes the asymptotic properties of the quantities
$z_k$ corresponding to a random vector on the moment space $
\mathcal{M}_{n}([0,\infty))$ with density $g_n^{(\gamma,\delta)}$
defined in (\ref{momdichte0infty}).
It is a well-known consequence of the asymptotic behavior of the
density of the Gamma distribution.
%
\begin{theorem} \label{asymzeta}
Suppose $\vec{m}_n $ is a random vector of moments on the moment space
$\mathcal{M}_n ([0,\infty))$ with density $ g_n^{(\gamma,\delta)}$,
where the $\gamma_k$ are fixed, $\delta_1=\cdots= \delta_n =n$, and let
$z_k^{(n)}$ denote the $k$th component of the vector $\vec{z}_n =
(z_1^{(n)},\ldots,z_n^{(n)})$. Then\vadjust{\goodbreak} the standardized
random variable $z_k^{(n)}$ is
asymptotically normal distributed, that is,
\[
\sqrt{n} \bigl(z_k^{(n)} - 1\bigr) \xrightarrow{n \rightarrow\infty
}^{\mathcal
{D}} \mathcal{N} (0,1).
\]
\end{theorem}

By Theorem~\ref{asymzeta} the vector
$\sqrt{n} (\vec{z}_k^{(n)}\! -\! \vec{1})\!=\! \sqrt{n} ((z_1^{(n)},\ldots,
z_k^{(n)})^T \!-\! (1,\ldots,1)^T)$
is asymptotically multivariate normal distributed. In
order to derive a corresponding statement of the random vector $\vec
{m}_k^{(n)} = \psi^{-1}_k(\vec{z}_k^{(n)} )$ we will use the Delta
method and study first the image of the vector
$(1,\ldots,1)^T$ under the mapping $\psi^{-1}_k$.
\begin{lem} \label{catalan}
Let $(c_n)_{n\geq1}$ denote the sequence of Catalan numbers, then
\[
\psi_n (c_1,\ldots,c_n ) = (1,\ldots,1)^T.
\]
\end{lem}
\begin{pf}
The proof presented here relies on the combinatorical
interpretation of the Catalan numbers and a recursive algorithm given in
\citet{skibinsky1968} to calculate the moments in terms of the variables
$z_k$. The $k$th Catalan number counts the paths in
$\en\times\en$ starting in $(0,0)$ and ending in $(2k,0)$, where one
is only allowed to make steps in the direction
$(1,1)$ or $(1,-1)$.
\citet{skibinsky1968} defines the triangular array $\{ g_{i,j}\}
_{i,j\geq0}$ by $g_{i,j}=0$ for $i>j$, $g_{0,j}=1$ and the recursion
%
%
\begin{equation}\label{skib}
g_{i,j} = g_{i,j-1} + z_{j-i+1} g_{i-1,j},\qquad 1\leq i \leq j.
\end{equation}
He showed that $g_{k,k} = m_k$. Consequently, if $z_i =1$ $(i=1,2,\ldots
)$, the quantity
$g_{k,k}$ is the number of paths through the lattice $\{ (i,j)\}
_{i,j\geq0}$, starting
in $(k,k)$ and ending in $(0,0)$, where in each vertex we can only
make steps upward
or to the left and where we are not allowed to cross the diagonal $\{
(i,i)\}$.
This number is exactly the $k$th Catalan number $c_k$.
\end{pf}
\begin{theorem} \label{konv0infty}
If the vector of random moments $\vec{m}_n \in\mathcal
{M}_n([0,\infty
))$ is governed by a law with density $g_n^{(\gamma,\delta)}$, where
$\delta_1=\cdots=\delta_n=n$, and the $\gamma_k$ are fixed, then the
projection $\vec{m}_k^{(n)}= \Pi_k^n (\vec{m}_n)$ of $\vec{m}_n$
onto the
first $k$ coordinates satisfies
\[
\sqrt{n} C^{-1} \bigl(\vec{m}_k^{(n)} - \vec{m}_k(\eta)\bigr) \xrightarrow{n
\rightarrow\infty}^{\mathcal{D}} \mathcal{N}_k(0,I_k),
\]
where the vector $\vec{m}_k(\eta) = (c_1,\ldots,c_k)^T$ contains the
first $k$ moments of the Marchenko--Pastur distribution, and $C$ is a
lower triangular matrix with entries $c_{1,1}=\cdots=c_{k,k}=1$, and
\[
c_{i,j} = \pmatrix{2i\cr i-j} - \pmatrix{2i\cr i-j-1},\qquad j < i .
\]
\end{theorem}
\begin{pf}
It suffices to calculate the Jacobi matrix
$
C = \frac{ \partial\psi_{k}^{-1}}{\partial\vec{z}_k} (\vec{z}_k^0)
$
of the mapping $\psi_{k}^{-1}$ at $\vec{z}_k^0 = (1,\ldots,1)^T$; then
the independence of\vadjust{\goodbreak} the recursion variables $z_k$ and Theorem
\ref{asymzeta} yield
\[
\sqrt{n} \bigl(\vec{m}_k^{(n)} - \vec{m}_k(\eta)\bigr) = \sqrt{n} C \bigl( \vec
{z}^{(n)}_k - \vec{z}_k^0\bigr) + o_P(1) \xrightarrow{n \rightarrow\infty
}^{\mathcal{D}}
\mathcal{N}_k(0,CC^T).
\]
Note that the moment $m_i$ depends only on $z_1,\ldots,z_i$, and
consequently $C$ is a lower triangular matrix.
To identify the entries of the matrix $C$ we consider the
triangular array $\{ g_{i,j}\} _{i,j\geq0}$ defined in (\ref{skib}).
For a fixed $r$ with $1\leq r \leq k$, we introduce the notation
$
u_{i,j} = \frac{ \partial g_{i,j}}{\partial z_r} (\vec{z}_k^0)
$
and obtain\vspace*{1pt} a new triangular array $\{ u_{i,j}\} _{i,j\geq0}$.
Obviously we have $u_{i,j}=0$ for $i>j$, and the other values of
$u_{i,j}$ are determined by
the initial condition $u_{0,j}=0$ and the recursion
%
%
\begin{equation} \label{reku}
u_{i,j} = u_{i,j-1} + u_{i-1,j} + \delta_{r,j-i+1} g_{i-1,j}^0,\qquad
1\leq i \leq j,
\end{equation}
were $\delta_{i,j}$ denotes the Kronecker symbol, and $g^0_{i,j}$ is
the coefficient in the recursion
(\ref{reku}), if all $z_i$ are equal to $1$, that is,
\[
g^0_{i,j} = \pmatrix{i+j\cr i} - \pmatrix{i+j\cr i-1}.
\]
The numbers $g^0_{i,j}$ are sometimes called generalized Catalan
numbers; see \citet{finu1976}.
An induction argument shows that the 
entries in the new triangular array are given by
%
%
\begin{equation} \label{loesungreku}
u_{i,j} =
\cases{\pmatrix{i+j\cr i-1} - \pmatrix{i+j\cr i-r-1}, &\quad if
$j-i\geq r$, \vspace*{2pt}\cr
\pmatrix{i+j\cr j-r} - \pmatrix{i+j\cr i-r-1}, &\quad if
$0 \leq j-i < r$.}
\end{equation}
With this identity we obtain for the entries of the matrix $C$
\[
c_{i,r} = \frac{ \partial m_i}{\partial z_r} (\vec{z}_k^0) = u_{i,i} =
\pmatrix{2i\cr i-r} - \pmatrix{2i\cr i-r-1}
\]
for $1\leq r \leq i$, which proves the assertion of the theorem.
\end{pf}

By the same arguments as in the compact case, the general class of
densities on the unbounded moment space $\mathcal{M}_n([0,\infty))$ for
which the result of Theorem~\ref{konv0infty} holds, is $g_n^{(\gamma
,\delta)}(\vec{m}_n) \prod_{k=1}^n \tilde{g}_k(z_k(\vec{m}_n))$. Here,
the functions $\tilde{g}_k$ have to be bounded and continuous and
positive at 1. An analogous result holds in the remaining case of
measures on $\er$.

We now consider the moment space $\mathcal{M}_{2n-1} (\er)$ and recall
the bijective mapping
(\ref{defxi}) from the interior of the moment space $\mathcal{M}_{2n-1}
(\er)$ onto $( \er\times\er^+ )^{n-1} \times\er$
corresponding to the range for coefficients in the recursive relation
of the orthogonal polynomials (\ref{polrek}).
The following results give the weak asymptotics of random recursion
coefficients and moments and correspond to Theorem~\ref{asymzeta} and
\ref{konv0infty}. The proof is omitted.
\begin{theorem} \label{asymanbn}
Let the random vector $\vec{m}_{2n-1} \in\mathcal{M}_{2n-1}(\er)$ be
governed by a law with density $h_{2n-1}^{(\gamma,\delta)}$
where $\gamma_k > -1$ is fixed and\vadjust{\goodbreak} $\delta_k=n$
$(k=1,\ldots,n)$. For
fixed $k$ denote by $b_k^{(n)}$ and by $a_k^{(n)}$ the $(2k-1)$th
component of the vector $\xi_{2n-1}(\vec{m}_{2n-1})$ and the $2k$th
component, respectively. Then
\[
\sqrt{2n} b_k^{(n)} \sim\mathcal{N} (0,1),\qquad \sqrt{2n}\bigl(a_k^{(n)} -1\bigr)
\xrightarrow{n \rightarrow\infty}^{\mathcal{D}} \mathcal{N} (0,1).
\]
\end{theorem}
\begin{theorem} \label{asymer}
Let the vector of random moments $\vec{m}_{2n-1} \in\mathcal
{M}_{2n-1}(\er)$ be governed by a law with density $h_{2n-1}^{(\gamma
,\delta)}$ where $\delta_k =
n$ $(k=1,\ldots,n)$. For $k \in\en$ denote by $\vec{m}_k^{(n)} = \Pi
_k^n(\vec{m}_{2n-1})$ the projection onto the first $k$ coordinates and
by $\vec{m}_k(\rho)= \Pi_k(0,c_1,0,c_2,\ldots
)$ the vector of the first $k$ moments of the semicircle law defined in
(\ref{marchpast}). Then
\[
\sqrt{2n} D^{-1} \bigl(\vec{m}_k^{(n)} - \vec{m}_k(\rho)\bigr) \xrightarrow{n
\rightarrow\infty}^{\mathcal{D}} \mathcal{N}_k(0,I_k),
\]
where $D$ is a $k\times k$ lower triangular matrix
with $d_{i,j}=0$ if $i+j$ is odd, and the remaining entries are given by
\[
d_{i,j} = \pmatrix{i\vspace*{2pt}\cr\dfrac{i-j}{2}} -
\pmatrix{i\vspace*{2pt}\cr\dfrac {i-j}{2} - 1}.
\]
\end{theorem}
\begin{rem}
It follows from Remark~\ref{ranspect}
that the moment density of the Jacobi ensemble is the moment density
investigated asymptotically in this section. Although for the random
matrix ensembles the parameters $\gamma_k,\delta_k$ depend on $n$, only
minor changes are necessary to obtain a weak convergence result for the
first $k$ moments. Note that the canonical moment $p_k^{(n)}$ follows a
Beta distribution with parameters behaving like $\frac{\beta}{2}n$. Therefore
$
\sqrt{4\beta n} (p_k^{(n)} -\frac{1}{2} ) \xrightarrow{n
\rightarrow
\infty}^{\mathcal{D}} \mathcal{N}(0,1),
$
and we obtain easily for the vector $\vec{m}_k(\mu_n)$ of the first $k$
moments of $\mu_n$ of the spectral measure of the Jacobi ensemble
defined in (\ref{jacew})
\[
\sqrt{4\beta n} A^{-1} \bigl(\vec{m}_k(\mu_n) - \vec{m}_k(\nu)\bigr)
\xrightarrow
{n \rightarrow\infty}^{\mathcal{D}} \mathcal{N}_k(0,I_k),
\]
where $\vec{m}_k(\nu)$ is the moment vector of the arcsine measure and
$A$ is the $k\times k$ matrix in (\ref{chakemsturesult}).

In particular, the moment convergence implies the weak convergence of
the spectral measure to the arcsine measure. This is also a consequence
of the well-known convergence of the empirical eigenvalue distribution
to the arcsine measure, since the (unscaled) moments of the spectral
measure have the same asymptotic behavior as the moments of the
empirical eigenvalue distribution. Therefore the fluctuations around
this limit in terms of the moments are Gaussian. The corresponding
results for the Laguerre and Hermite ensemble are omitted for the
sake of brevity.
\end{rem}
\begin{rem}
Already in the compact case there are interesting results regarding a
functional limit theorem. In particular,\vadjust{\goodbreak} \citet{detgam2007} proved the
convergence of the standardized range process $(m^+_{\lfloor nt\rfloor}
- m^-_{\lfloor nt\rfloor})_t$\vspace*{2pt} to a functional of the Brownian motion in
the Skorohod topology. With the distributions on the unbounded moment
spaces, the question arises whether a corresponding result exist in
these cases. Interesting processes are, for example, the moment
difference $(m_{\lfloor nt\rfloor} - m^-_{\lfloor nt\rfloor})_t$ for
measures on $[0,\infty)$ and for measures on the whole real line the
integrals over orthogonal polynomials as in formulas (\ref{wall1}) and
(\ref{wall2}). 
For the sake of brevity, we do not discuss functional limit theorems in
the unbounded cases and defer these interesting questions to our future
research.
\end{rem}

\section*{Acknowledgments}

The authors are grateful to Martina Stein who typed parts of this
paper with considerable technical expertise. The authors would also
like to thank two anonymous referees for their constructive comments on
an earlier version of this paper.


%

\printaddresses

\end{document}